%% file: 2-cable.tex
\documentclass[11pt]{amsart}
\font\emailfont=cmtt10

\usepackage{amsmath,amsthm,amsfonts,amscd,flafter,epsf}
\usepackage{graphics}
\usepackage{epsfig}
\usepackage{psfrag}

\input{mac}

\commentable{

\title[{Khovanov homology of the $2$-cable detects the unknot}] 
{Khovanov homology of the $2$-cable detects the unknot}

}

\author[Matthew Hedden]{Matthew Hedden}
\address{Department of
Mathematics, Massachusetts Institute of Technology, MA \newline
\indent{\emailfont{mhedden@math.mit.edu}}}

\thanks{The author was partially supported by NSF Grant DMS-0706979}

\begin{document}

\begin{abstract}
We prove that the Khovanov homology of the $2$-cable detects the unknot. A corollary is that Khovanov's categorification of the $2$-colored Jones polynomial detects the unknot.  \end{abstract}

\maketitle
\section{Introduction}

A central problem in the study of quantum knot invariants is to determine if the Jones polynomial detects the unknot \cite{Kirby}. Even if it does not, one could hope that the Jones polynomial of {\em cables} of a knot could detect the unknot.    To make this precise, recall that the $n$-cable of $K$, denoted $K^n$, is the $n$-component link obtained by taking $n$ parallel copies of $K$ (this  operation depends on a framing, which we will always assume to come from a Seifert surface).  Since isotopic knots have isotopic cables and the $n$-cable of the unknot, $U$, is the trivial $n$-component link, $U^n$, one could try to show that if $K$ is non-trivial then $V_{K^n}\ne V_{U^n}$ (here, $V$ denotes the Jones polynomial).


In \cite{Khovanov}, Khovanov introduced a link invariant which takes the form of a bi-graded homology theory.  The graded Euler characteristic of his invariant is the Jones polynomial. Thus Khovanov's invariant contains as much information as the Jones polynomial, and one can relax the above questions to the context of Khovanov homology.  The purpose of this note is to show that the Khovanov homology of the $2$-cable detects the unknot.  

\begin{theorem}\label{thm:main}  Let $K$ be a knot, and let $K^2$ denote its $2$-cable. Then 
$ \mathrm{rk} \ {{\mathrm{Kh}}} (K^2)=4 $
 if and only if $K$ is the unknot. Here, ${\mathrm{Kh}}$ denotes the Khovanov homology (with coefficients in $\Z/2\Z$).  
\end{theorem}
The key idea behind this theorem is that there is a spectral sequence from (reduced) Khovanov homology of $K^2$ to the \os \ Floer homology of the branched double cover of $K^2$. The Thurston norm of this latter manifold is non-trivial in the case when $K$ is non-trivial.   Since \os \ Floer homology detects the Thurston norm, this forces the rank of Floer homology, and hence the rank of Khovanov homology, to grow.  

We should remark that the Khovanov homology of cables is closely related to the so-called ``colored" Khovanov invariants, whose graded Euler characteristics are the colored Jones polynomials (see \cite{Khovanov2} for definitions).  The following is an immediate corollary of our proof.

 \begin{cor}\label{cor:color} Khovanov's categorification of the $2$-colored Jones polynomial detects the unknot.  
\end{cor}

The present work should be compared to that of Grigsby and  Wehrli \cite{Grigsby}, who show that there are spectral sequences from Khovanov's categorification of the {\em reduced} colored Jones polynomial to an invariant coming from Juhasz's sutured Floer homology \cite{Juhasz}.  Though both results use spectral sequences from a Khovanov style invariant to a Floer invariant, the actual techniques are quite different.  

It would also be interesting to compare our result with that of Anderson \cite{Anderson}, who has announced a proof that the unknot is detected by the collection of the Jones polynomials of {\em all} cables of a knot.

\smallskip
\smallskip
\noindent{\bf{Acknowledgment:}} The inspiration for this note came from a talk I saw given by Stephan Wehrli at the ``Knots in Washington XXVI" conference at George Washington University and from a conversation with Eli Grigsby at the 2008 AMS sectional meeting at LSU.  I am indebted to them for their willingness to discuss their results. I also thank Liam Watson for interesting conversations, which led to \cite{Tangle}.

\section{Proof of Theorem}
\label{sec:proof} 
\noindent All homology theories will have coefficients in $\Z/2\Z$.
\bigskip

\noindent {\bf Proof of Theorem \ref{thm:main} }
Assume that $K$ is a non-trivial knot.

 Let $\Sigma(L)$ denote the double cover of $S^3$, branched along a link $L$. The proof begins by recalling Theorem $1.1$ of \cite{Branched}, which says that there exists a spectral sequence whose $E_2$ term 
consists of Khovanov's reduced homology of the mirror of $L$,  and which converges to the \os \ ``hat" Floer homology of  $\Sigma(L)$.  An immediate corollary (Corollary $1.2$ of \cite{Branched}) is that $$\mathrm{rk}\ \widetilde{\mathrm{Kh}}(L) \ge \mathrm{rk} \ \HFa(\Sigma(L)).$$   Here, $\widetilde{\mathrm{Kh}}$ denotes the reduced Khovanov homology and $\HFa$ the ``hat"  Floer homology (see \cite{Khovanov} resp. \cite{HolDisk} for relevant definitions). 

Thus it remains is to bound the rank of $\HFa(\Sigma(K^2))$. To achieve this, we first observe that $$\Sigma(K^2)\cong S^3_0(K\#K),$$ where $S^3_0(K\#K)$ denotes the manifold obtained by zero-framed surgery on the connected sum of $K$ with itself \footnote{Strictly speaking, we should reverse the string orientation of one of the copies of $K$ when taking the connected sum. Since our argument is insensitive to this discrepancy we suppress it from the notation.}.  This can be seen, for instance,  from the technique of Akbulut and Kirby \cite{Kirby} for finding a Kirby calculus description of the branched double cover of $L$ in terms of a Seifert surface.  Specifically, apply the algorithm indicated by Figure $4$ of \cite{Kirby} to the annulus which $K^2$ bounds.  

We now claim that if $K$ is non-trivial, then rk $\HFa(S^3_0(K\#K))\ge 6$.  Essentially, this follows from the fact that the Thurston norm of this manifold is non-trivial, together with the fact that Floer homology detects the Thurston norm.

Let us be more precise.  First, we have the elementary fact that the genus of $K\#K$ is $2g$, where $g$ is the genus of $K$.  Theorem $1.2$ of \cite{GenusBounds} says that knot Floer homology detects the genus of $K$. In particular, it implies that $\HFKa(K\#K,2g)\ne 0$. 

Now let $\HFp(S^3_0(K\# K),i)$ denote the ``plus" Floer homology group associated to the $\SpinC$ structure $\spinc_{i}\in\SpinC(S^3_0(K\# K))$ whose Chern class satisfies $ \langle c_1(\spinc_{i}),[\widehat{F}]\rangle = 2i.$  Here, $[\widehat{F}]\in H_2(S^3_0(K\# K);\Z)\cong\Z$ is the generator obtained by capping off the oriented Seifert surface of $K\#K$ with the meridional disk of the surgery torus. 

Since $2g>1$, we can apply Corollary $4.5$ of \cite{Knots} to show that $$\HFp(S^3_0(K\# K),{2g-1})\cong \HFKa(K\#K,2g).$$
 
An application of the  long exact sequence relating $\HFp$ to $\HFa$ (Lemma $4.4$ of \cite{HolDisk}) then shows that  $\HFa(S^3_0(K\# K),2g-1)\ne 0$.  Since the Euler characteristic of this group is zero (Proposition $5.1$ of \cite{HolDiskTwo}), we have $$\mathrm{rk}\ \HFa(S^3_0(K\# K),2g-1)\ge 2.$$

Floer homology enjoys a symmetry relation under the action of conjugation on $\SpinC$ structures (Theorem $2.4$ of \cite{HolDiskTwo}).  In the present context, this implies $$\mathrm{rk}\ \HFa(S^3_0(K\# K),-2g+1)\ge 2.$$

Finally, Theorem $10.1$ of \cite{HolDiskTwo}, together with the long exact sequence relating $\HFp$, $\HFm$, and $\HFinf$  shows that $\HFp(S^3_0(K\# K),0)\ne0$.  By the same reasoning as above, we obtain $$\mathrm{rk}\ \HFa(S^3_0(K\# K),0)\ge 2,$$
\noindent and have proved the claim.

Since the rank of the reduced Khovanov homology of the $2$-component unlink is $2$, we have shown that the reduced Khovanov homology of the two-cable detects the unknot.  To obtain the statement in terms of the (unreduced) Khovanov homology, we observe that (since we are working with $\Z/2\Z$ coefficients), 
$$ \mathrm{Kh}(K) \cong \widetilde{\mathrm{Kh}}(K)\otimes V,$$
where $V$ is the rank $2$ vector space obtained as the unreduced Khovanov homology of the unknot.   This follows, for instance, from \cite{Shumakovitch}.   Thus, the Khovanov homology of the $2$-cable of a non-trivial knot has rank at least $12$.   \ \ \ \ \ \ \ \ \ \ \ \ \ \ \ \ \ \ \ \ \ \ \ \ \ \ \ \ \ \ \ \ \ \ \ \ \ \ \ \ \ \ \ \ \ \ \ \ \ \ \ \  \ \ \ \ \ \ \ \ \ \ \ \ \ \ \ \ \ \ $\square$

\begin{remark} The above argument extends in a straightforward manner to show that the Khovanov homology of the $2n$-cable detects the unknot.  A slight modification should handle the case of odd cables. In this case one must show that the Thurston norm of a certain manifold obtained by surgery on $\Sigma(K)$ is non-trivial. 
\end{remark} 

\noindent{\bf Proof of Corollary \ref{cor:color}:}
 From its definition \cite{Khovanov2}, the rank of Khovanov's categorification of the $2$-colored Jones polynomial differs from that of the Khovanov homology of the $2$-cable by at most $1$.  Thus the rank of Khovanov's categorification of the $2$-colored Jones polynomial of a non-trivial knot is at least $11$, whereas it is $3$ for the unknot.   \ \ \ \ \ \ \ \ \ \ \ \ \ \ \ \  \ \ \ \ \  \ \ \ \ \ \ \ \ \ \ \ \ \ \ \  $\square$

\end{document}

%% file: mac.tex
\newcommand\commentable[1]{#1}

\newcommand{\HF}{HF}

\newtheorem{theorem}{Theorem}[section]

\newtheorem{cor}[theorem]{Corollary}

\newtheorem{remark}[theorem]{Remark}

\def\endproofof{\relax\ifmmode\expandafter\endproofmath\else
  \unskip\nobreak\hfil\penalty50\hskip.75em\hbox{}\nobreak\hfil\bull
  {\parfillskip=0pt \finalhyphendemerits=0 \bigbreak}\fi}
\def\endproofofmath$${\eqno\bull$$\bigbreak}

\def\endproof{\relax\ifmmode\expandafter\endproofmath\else
  \unskip\nobreak\hfil\penalty50\hskip.75em\hbox{}\nobreak\hfil\bull
  {\parfillskip=0pt \finalhyphendemerits=0 \bigbreak}\fi}
\def\endproofmath$${\eqno\bull$$\bigbreak}
\def\bull{\vbox{\hrule\hbox{\vrule\kern3pt\vbox{\kern6pt}\kern3pt\vrule}\hrule}}

\newcommand{\Z}{\mathbb{Z}}

\newcommand\SpinC{\mathrm{Spin}^c}

\newcommand\relspinc{\underline{\spinc}}

\newcommand\ModSphere{\ModFlow\left({\mathbb S}\longrightarrow 
\Sym^{g-1}(\Sigma_{1})\times \Sym^2(\Sigma_{2})\right)}
\newcommand\ModSpheres\ModSphere

\newcommand\HFp{\HFb}

\newcommand\HFm{\HF^-}

\newcommand\HFinf{HF^\infty}

\newcommand\HFa{\widehat{HF}}
\newcommand\HFb{HF^+}

\newcommand\UnparModSp{\widehat \ModSp}
\newcommand\UnparModFlow\UnparModSp
\newcommand\Mod\ModSp

\newcommand{\spinc}{\mathfrak s}

\newcommand\ModMaps{\mathcal M}
\newcommand\ModSp\ModMaps

\newcommand\spincrel\relspinc

\newcommand\HFK{HFK}

\newcommand\HFKa{\widehat\HFK}

\newcommand\Dual{\mathcal D}
\newcommand\Duality\Dual

\newcommand\os{{Ozsv{\'a}th-Szab{\'o}}}